\documentclass[11pt,batang]{article}

\usepackage{amsmath}
\usepackage{amsfonts}
\usepackage{amssymb}
\usepackage{amsthm}
\usepackage{graphicx}
\usepackage{hyperref}

\usepackage{color}
\usepackage{enumitem}
\usepackage{url}
\usepackage[margin=4cm]{geometry}
\usepackage{epigraph}
\theoremstyle{plain}

\newtheorem{theorem}{Theorem}[section]

\newtheorem{definition}[theorem]{Definition}

\numberwithin{equation}{section}

\begin{document}

\title{
 \begin{flushright}
\includegraphics[width=8 cm]{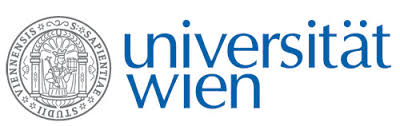}
\end{flushright}\vskip 2cm
 \begin{Huge} Habilitationsschrift\end{Huge}\vskip 2cm \LARGE Affine Processes\vskip 3cm \large Dipl.-Ing. Dr.rer.nat. Eberhard Mayerhofer\vskip 4 cm
\begin{flushleft} \hskip 1cm Dublin, April 2014\end{flushleft}
}
\date{}

\maketitle

\tableofcontents
\newpage

\section{Summary}

 We put forward a complete theory on moment explosion for fairly general state-spaces [AP1], [AP2] and [AP9]. This includes a characterization of the validity of the affine transform formula in terms of minimal solutions of a system of generalized Riccati differential equations.

Also, we characterize the class of positive semidefinite processes [AP4], [AP5], and provide existence of weak and strong solutions for Wishart SDEs. As an application of [AP4],  we answer a conjecture of M.L. Eaton on the maximal parameter domain of non-central Wishart distributions [AP6].

The last chapter of this thesis comprises three individual works on affine models, such as a characterization of the martingale property of exponentially affine processes [AP3], an investigation of
the jump-behaviour of processes on positive semidefinite cones [AP7] and an existence result for transition densities of multivariate affine jump-diffusions and their  approximation theory in weighted Hilbert spaces [AP8].
\newpage
\section{Publications}
\begin{itemize}
\item [] {\bf [AP9]} (with Martin Keller-Ressel) {\it On exponential moments of affine processes}, to appear in Annals of Applied Probability (2014).
\item [] {\bf [AP8] }(with D. Filipovi\'c and P. Schneider) {\it Density approximations for multivariate affine jump-diffusion
processes},  Journal of Econometrics, \\ Vol. 176 (2013), 93--111.
\item [] {\bf [AP7] }{\it Affine processes on positive semidefinite $d \times d$ matrices have jumps of finite variation in dimension $d>1$}, Stochastic Processes and Their Applications 122 (2012), Issue 10, 3445--3459.
\item [] {\bf [AP6] }{\it On the existence of non-central Wishart
distributions}, Journal of Multivariate Analysis 114 (2013), 448--456.
\item [] {\bf  [AP5] }(with Oliver Pfaffel and Robert Stelzer), {\it On strong solutions of positive definite jump-diffusions}, Stochastic Processes and their Applications 121 (2011), No. 9, 2072--2086.
\item [] {\bf  [AP4] }(with C. Cuchiero., D. Filipovi\'c and J. Teichmann), {\it Affine processes on
positive semidefinite matrices}, Annals of Applied Probability,
Vol. 21, No. 2, 397--463 (2011).
\item [] {\bf  [AP3] }(with Johannes Muhle-Karbe and Alexander G. Smirnov), {\it A characterization of the martingale property of
exponentially affine processes}, Stochastic Processes and their Applications,  121 (2011), No. 3, 568--582.
\item [] {\bf  [AP2]}(with M. Keller-Ressel and A. G. Smirnov),
{\it On convexity of solutions of ordinary differential equations,}
Journal of Mathematical Analysis and Applications 368 (2010), 247--253.
\item [] {\bf [AP1]} (with D. Filipovi\'c), {\it Affine diffusion processes:
theory and applications}, Advanced Financial Modelling, Radon  Series
of Computational and Applied Mathematics, Vol. 8, 125--164, Walter de Gruyter, Berlin, 2009.
\end{itemize}
\newpage
\section{Preface}

Most of this work resulted from collaborations with colleagues who
were in Vienna.

I thank Michael Kunzinger, my doctoral supervisor at the University of Vienna, for leading me to my current mathematical standards during the project on ``The wave equation on singular space-times". I am also grateful to Damir Filipovi\'c and Josef Teichmann, who are responsible for my directional change in research after my PhD thesis. Josef I admire for his research oriented lectures in Stochastic Analysis at TU Vienna, and Damir I owe for leading me very efficiently into the field of Affine Processes and sharing open research problems with me. I have very much enjoyed working closely with Christa Cuchiero, Martin Keller-Ressel and Paul G. Schneider.

This habilitation thesis is submitted during my lecturership in financial mathematics at Dublin City University. I thank Paolo Guasoni for the fine academic environment in Dublin and the financial support through ERC Grant 278295.


\newpage

\section{Introduction}

In the vicinity of finance, many questions of a mathematical nature arise. This is due to the fact that the mainstream of financial models is probabilistic and even only remotely realistic stochastic models bring complex mathematical problems. A sound theoretic understanding of these provides the basis for reliable applications.

Affine Markov processes comprise most of the continuous-time stochastic models in Finance. This work deals with different theoretical aspects, in particular with existence of certain affine Markovian semigroups  and their analytic and stochastic properties. The key property of these processes is a particular simple action of the affine semigroups
on exponentially affine functions, which makes them ``analytically tractable'' to a certain extent\footnote{ Analytic tractability often expresses
the fact that one can derive quantities in closed or semi-closed form, thus avoiding time-extensive simulations or complex approximations.}:

\begin{definition}
A stochastically continuous Markov process $(X,(\mathbb P^x)_{x\in D})$ with state-space $D\subset \mathbb R^d$ is affine, if
there exist functions $\phi: \mathbb R_+\times i\mathbb R^d\rightarrow \mathbb C$
and $\psi: \mathbb R_+\times i\mathbb R^d\rightarrow \mathbb C^d$ such that for all $x\in D$, $t\in\mathbb R$ and $u\in i\mathbb R^d$ we have\footnote{Angle brackets denote
an inner product on $\mathbb R^d$.}
\begin{equation}\label{eq: affine property}
\mathbb E[e^{\langle u, X_t\rangle}\mid X_0=x]=e^{\phi(t,u)+\langle \psi(t,u),x\rangle}.
\end{equation}
\end{definition}
Many well-known classes of stochastic processes satisfy \eqref{eq: affine property}, e.g.,
\begin{itemize}
\item L\'evy processes, $L_t$ for which $\psi(t,u)\equiv t u$, and L\'evy driven OU-type processes  on $\mathbb R^d$ or subsets thereof, with linear drift $\beta\in\mathbb R^{d\times d}$, i.e.
\[
X_t=x+\int_0^t\beta^\top X_s ds+L_t.
\]
In this case we have,
\[
\psi(t,u)=e^{\beta t}u.
\]
\item Square Bessel processes, first introduced by William Feller \cite{feller51}, which solve the stochastic integral equation
\begin{equation}\label{eq: cir}
X_t=x+\int_0^t (b+\beta X_s)ds+\sigma\int_0^t \sqrt{X_s}dB_s, \quad x\geq 0,
\end{equation}
where $B$ denotes a one-dimensional standard Brownian motion. In the finance literature, this model is  used by Cox, Ingersoll and Ross \cite{cir85}. Among the mentioned examples, this is the first
one with a non-trivial function $\psi$, satisfying
\[
\partial_t\psi(t,u)=\frac{\sigma^2}{2}\psi^2(t,u)+\beta\psi(t,u),\quad \psi(0,u)=u.
\]
\item Stochastic volatility models, such as Bates' \cite{Bates2000}, Heston's \cite{Heston1993} and Barndorff-Nielsen \& Shepard's \cite{Barndorff2001} (state space $D=\mathbb R_+\times \mathbb R$).
\item Multi-variate interest rate models such as Dai-Singleton class \cite{daisingleton00}, \cite{duffiekan96}, or the jump-diffusion models of Duffie, Pan and Singleton \cite{duffie2000transform} (state space $D=\mathbb R_+^m\times \mathbb R^n$).
\item Branching processes with immigration, \cite{kw71} (state space $D=\mathbb R_+^m$. All the previously mentioned classes of affine processes are covered by the general theory
of Duffie, Filipovi\'c and Schachermayer \cite{dfs}.
\item Most relevant for this work are the Wishart processes introduced by Bru \cite{bru}, which are matrix generalizations of \eqref{eq: cir}, with range on the cone of positive semidefinite matrices $D=S_d^+$ (pure diffusions), and the matrix L\'evy-subordinator driven, Ornstein-Uhlenberck-type processes introduced by Barndorff-Nielsen and Stelzer, cf.~\cite{bierundstelze}. Using any of these as a covariance process yields a quite straightforward generalization of Heston's, or of the BNS model. The sophisticated issue in the latter case is the choice of non-trivial dependence structure with one-dimensional marginal distributions given by the BNS model. Useful suggestions for modeling such dependencies have been put forward by Perez-Abreu and Stelzer \cite{perez2012class}.

\end{itemize}

The main aim of the present thesis is to {\it design methods} and provide theories {\it which are inherently affine:}  They
are tailored to affine processes, in that they try to use, more than anything else, two of their key properties:
\begin{enumerate} \item [(a)]\label{key prop 1} The exponentially affine form
\eqref{eq: affine property} of the characteristic function. Stochastic continuity implies regularity of the characteristic exponents ($\phi,\psi$)  \cite{regu1}, \cite{regu2}, from which it may be concluded that for $u\in i\mathbb R^d$
the so-called generalized Riccati differential equations are satisfied
\begin{align}\label{eq ric1}
\partial_t\psi(t,u)&=R(\psi(t,u),\quad \psi(0,u)=u,\\\label{eq ric2}
\partial_t\phi(t,u)&=F(\psi(t,u),\quad \phi(0,u)=0.
\end{align}

 The right sides of equations \eqref{eq ric1}--\eqref{eq ric2} are component-wise of L\'evy-Khintchine form on $\mathbb R^d$, satisfying special parametric constraints which crucially depend on the geometry of the state space $D$.

\item [(b)]\label{key prop 2} The polynomial character of conservative processes \cite{CKRT}: The Markovian semigroup, whenever defined, maps the finite dimensional real vector space $P_k$ of polynomials of degree $\leq k$ into itself.\footnote{ This special property is shared by e.g., Jacobi Processes, which have a quadratic diffusion coefficient and thus are not affine.}  The polynomial property implies, that the Markovian semigroup action on $P_k$ reduces to that of a finite dimensional semigroup, one that can be expressed as a matrix exponential.
\end{enumerate}

Also, the {\it covered topics are} exclusively dedicated to the {\it multivariate} nature of affine processes. In fact, some of the theories go far beyond what might be considered as feasible to be understood in a general Markovian framework, e.g., our characterization [AP3], [AP9] of true exponentially affine martingales in terms of deterministic conditions, see section \ref{sec martingales}.

In the following I describe the three parts of the thesis in more detail, using minimal references to related work. The aim is to explain  intuitively rather than theoretically
my work and also to connect the papers among each other rather than to related research; detailed references to the literature are provided then
in the respective papers.

\subsection{Moment Explosion}\label{sec: Moment Explosion}
[AP1], [AP2] and [AP9] are dedicated to studying the moment generating function
\[
g(t,u,x):=\mathbb E[e^{\langle u, X_t\rangle}\mid X_0=x]
\]
for real initial data $u\in\mathbb R^d$, and also, for general complex data $u\in\mathbb C^d$. We focus on relating moment explosion\footnote{Moment explosion occurs when for fixed $t>0$ the MGF $u\mapsto g(t,u,x)$ ceases to be finite.} to the explosion
in the associated Riccati differential equations in a meaningful way. One might expect that if either side of \eqref{eq: affine property}
is finite, then also the other is, and equality holds. Furthermore $\phi,\psi$ solve the generalized Riccati differential equations
 \eqref{eq ric1}--\eqref{eq ric2}. We call this the validity of the affine transform formula (ATF for short). Problems that one encounters in achieving a result like this include the following:
\begin{enumerate}
\item [(a)] The generalized Riccati differential equations \eqref{eq ric1}--\eqref{eq ric2} need not be well posed for arbitrary initial data: The right sides $F,R$ are only continuous, but not Lipschitz, in general. Uniqueness issues are inevitable.
\item [(b)] Explosion to solutions of initial value problems (IVPs) is well understood in terms of the standard existence and uniqueness result for ODEs. However, explosion in the initial data is not generally understood, i.e. the behaviour $u\mapsto g(t,u,x)$.
\item [(c)] Can Markov processes recover from Moment explosion? Or can we conclude from $g(t,u,x)=\infty$ that also $g(t+h,u,x)=\infty$, each $h>0$. \footnote{Apart from general Markov chains in continuous time (and from affine processes), this question remains unsolved.}
\end{enumerate}

For pure diffusion models on state spaces of the type $D=\mathbb R_+^m\times\mathbb R^n$, most of the technicalities, in particular concerning (a)  and (c)\footnote{That diffusions do not recover from moment explosion is a mere consequence of [AP1] and not technically necessary for proving the main characterization in [AP1].}, are not an issue. Inspired by a  result of  Glasserman \& Kim \cite{KG2010} for a particular class of affine term structure models, Damir Filipovi\'c and I [AP1] prove the validity of the ATF. The characterization of real moments is solved by using multivariate ODE comparison results in conjunction with a certain blow-up-lemma, which elaborates on (b).

The jump-diffusion case has been solved jointly with Martin Keller-Ressel in [AP9] for general convex state spaces. Uniqueness of solutions of
\eqref{eq ric1}--\eqref{eq ric2} is given up and is replaced by a concept of ``minimal solutions''. The paper [AP2] prepared a related result in [AP9] on the existence of complex exponential moments and the validity of the affine transform formula, by showing that IVPs with convex, quasi-monotone increasing  (with respect to a proper, closed convex cone) right side, 
have convex dependence on initial data. The ATF for complex arguments requires sophisticated a-priori estimates for the generalized Riccati differential equations,
hence for this part we elaborate on the most relevant state spaces $D=\mathbb R_+^m\times\mathbb R^n$ and $S_d^+$.

An important application of [AP9] to exponentially affine martingales and their relation to conservative affine semigroups as in [AP3] is described in section \ref{sec martingales}.

\subsection{Positive Semidefinite Processes}
[AP4] and [AP5] comprise foundations for affine processes on $D=S_d^+$, the state-space of positive semidefinite matrices of arbitrary dimension $d\geq 1$. [AP4] characterizes affine Markovian semigroups in terms of a parametric description of their infinitesimal generators. One implication of this theory is the existence of weak solutions for a large class of associated stochastic differential equations. Strong solutions for strictly positive definite affine jump-diffusions are derived in [AP5]. The paper [AP8] then
applies [AP4] to a characterization of the maximal parameter domain of non-central Wishart distributions.

 The theory of [AP4] complements previous research on existence by Kawazu and Watanabe \cite{kw71} ($D=\mathbb R_+$) and
by Duffie, Filipovi\'c and Schachermayer \cite{dfs} (for canonical state-spaces $\mathbb R_+^m\times \mathbb R^n$). There is, however, a special feature of affine semigroups on
$S_d^+$ which distinguishes them from those on canonical state spaces. In general, their transition functions are not infinitely divisible. 
The lack of infinite divisibility made a new approach to existence necessary. Let me give a little more insight into this issue. The principal symbol of the infinitesimal generator is given by
\[
p(x,\xi)=\text{tr}(\xi \alpha \xi x),\quad \xi,x\in S_d^+,
\]
where $\alpha\in S_d^+$ is the so-called diffusion coefficient of $X$. This entails that on a possibly enlarged probability space,
 $X$ can be written as the weak solution to an SDE of the form
\begin{equation}\label{weakest eq}
dX_t=(b+B^\top(X_t))dt+\sqrt{X_t}dW_t Q+Q^\top dW_t^\top \sqrt{X_t} +dL_t,\quad X_0\in S_d^+,
\end{equation}
where $Q$ is a $d\times d$ matrix satisfying $Q^\top Q=\alpha$,
where $W$ is a standard $d\times d$ Brownian motion and $L$ is a pure jump-process with affine jump intensity, and $B$ is a linear drift, inward pointing at the boundary of $S_d^+$. We know that in the absence of jumps $L=0$ and
for $B=0, Q=I$, and $b=\delta I$, $2\delta>d-1$, $X$ is a process introduced by \cite{bru}. That, in turn, has a transition function which is non-central Wishart. However, Wishart laws are, in general,
not infinitely divisible (this is a classical finding of L\'evy \cite{levywishart} for a special case on $S_2^+$, and generally established by \cite{PR} and others). This example already
suggests a general lack of infinite divisibility of the transition functions of APs on $S_d^+$.

The Feller property implies that the constant drift $b$ must satisfy the condition
\begin{equation}\label{eq drift}
b-(d-1) \alpha\in S_d^+,
\end{equation}
hence it is not sufficient that $b\succeq 0$, unless $d=1$ or $\alpha=0$ (the pure jump-case).\footnote{Note that such a drift condition does not exist in the case of canonical state spaces; there we have $b\in D$, see \cite{dfs}.}
We show in ([AP4], Theorem 2.9) that infinite decomposability is equivalent to having $\alpha=0$ or $d=1$ (the reason being that $b/n$ would violate the drift condition \eqref{eq drift} for sufficiently large $n$). A particular implication is that only in these cases,  affine jump-diffusions can be approximated weakly by pure jump-processes, as is in the L\'evy case (because infinite divisible distributions are closed under weak convergence).

A consequence of the Markov theory established in [AP4] is the existence of weak solutions of stochastic differential equations of the form \eqref{weakest eq}, for $L$ being
a matrix-valued L\'evy subordinator. In [AP5] we consider similar SDEs and provide existence results for strong solutions, under a more stringent drift condition, namely $b-(d+1)\alpha\in S_d^+$. The key to strong existence is a result on boundary non-attainment: for initial data in $S_d^{++}$, the cone of positive definite matrices, a process does not hit the boundary in finite time, given the drift condition.  The knowledge concerning boundary non-attainment is actually very useful
for obtaining affine structure preserving, equivalent changes of measures, \cite{AP10}.

[AP8] uses the Feller property of affine processes to derive information on the maximal parameterization of non-central Wishart distributions. Wishart processes, are solutions to
\eqref{weakest eq} when $L=0$ and $B(X)=\beta X+X\beta^\top$, with a $d\times d$ matrix $\beta$, and $b=2p \alpha$ ($2p\geq d-1$). Their transition function
has Laplace transform
\[
\int_{S_d^+} e^{-\text{tr}(u \xi)}p_t(x,d\xi)=\det(I+\sigma_t^\beta(u))^{-p} e^{-\text{tr}(u(I+\sigma_t^\beta(\alpha)u)^{-1}\omega_t^\beta(x))}.
\]
Hence $p_t(x,d\xi)$ is Wishart distributed, with shape parameter $p$, non-centrality parameter
$\omega_t^\beta(x)\in S_d^+$ and scale parameter  $\sigma_t^\beta(\alpha)\in S_d^+$ (the detailed definitions are found in [AP8]). An indirect proof, which constructs
a whole Markov transition function from a single distribution, and uses the constraint  $2p\geq d-1$, lets us
conclude that for any Wishart distribution with parameters $p,\sigma,\omega$ we must have
\begin{equation}\label{eq Gindikin}
p\in \{0,1,\dots,(d-2)/2\}\cup [(d-1)/2,\infty)
\end{equation}
and
\begin{equation}\label{eq rank}
\text{rank}(\omega)\leq 2p+1
\end{equation}
The necessary condition \eqref{eq Gindikin} answers a long standing conjecture by M.L. Eaton in the positive, which has only been previously solved by special function methods  for the $\text{rank}(\omega)=1$ case \cite{PR}. The rank condition \eqref{eq rank} is new. We further conjectured that for $2p<d-1$ we must actually have  
\begin{equation}\label{eq rank1}
\text{rank}(\omega)\leq 2p.
\end{equation}
This conjecture has recently been confirmed by \cite{LetacMassam} with means inspired by ours, but not using any reference to stochastic analysis.

\subsection{Properties of Affine Processes}
 [AP3], [AP6] and [AP7] study specific properties of Affine Processes or functionals thereof.
\subsubsection{Exponentially Affine Martingales}\label{sec martingales}
Several applications in financial mathematics require the knowledge, under which conditions an exponentially affine functional
\[
g(X_t):=e^{\langle \theta, X_t\rangle},\quad t\in [0,T]
\]
for $\theta\in\mathbb R^d$ is a true martingale. Two papers in this thesis characterize the martingale property of $g(X_t)$:
\begin{itemize}
\item[] [AM3] characterizes moment explosion and gives as an application a link of the martingale property of stochastic exponentials
to non-explosion of a related affine process. 
\item [] [AM9] characterizes the validity of the ATF and gives as application a link of the martingale property of exponentially affine functionals to 
moment explosion. This result is for general convex state spaces. 
\end{itemize}
Both moment explosion and the validity of the ATF are characterized in terms of a unique minimal solution to a system of generalized Riccati differential equations. These results can
actually be linked, since stochastic and standard exponentials of affine processes are related (cf [AM3], Remark 4.5 (i)). 
We are not going to perform the tedious calculations, which require an extension of state-space. Nevertheless, I think it is illuminating to demonstrate the relation of
validity of the ATF and non-explosion of an associated AP,  for standard exponentials of affine processes with state-space $D=\mathbb R_+$:

Let $X$ be a conservative affine process on $\mathbb R_+$ with infinitesimal generator
\begin{align*}
\mathcal Af(x)&=\alpha x f''(x)+(b+\beta x)f'(x)\\&+\int_{D\setminus \{0\}}(f(x+\xi)-f(x)-f'(x)\chi(\xi)) (m(d\xi)+x\mu(d\xi))
\end{align*}
where $\chi$ is a truncation function. Assume for $\theta\in\mathbb R$ and $x>0$ we have
\[
g(X_T):=\mathbb E[e^{\theta X_T}\mid X_0=x]<\infty.
\]
By  ([AP9], Theorem 2.14) we infer the existence of unique minimal solutions of the Riccati differential equations 
\eqref{eq ric1}--\eqref{eq ric2} on $[0,T]$ with initial data $\phi(0)=0,\psi(0)=\theta$. Furthermore, if $M_t:=e^{\theta X_t}$ is a $\mathbb P^x$ martingale, then $t\mapsto g(X_t)$ has constant expectation, hence we must have
that $\phi(t,\cdot)\equiv 0,\psi(t,\cdot)\equiv \theta$ and $R(\theta)=F(\theta)=0$ (see also [AP9], Remark 3.2 (2)). We may introduce by exponential tilting an affine Markov process  $(X,\mathbb Q_x)_{x\in D}$ with transition function\footnote{ For a complete proof in a more general
setting, in particular of the Markov property of $(X,\mathbb Q_x)_{x\in D}$, see Theorem 4.14 in \cite{keller}.}
\begin{equation}\label{tada}
q_t(x,d\xi):=\frac{e^{\theta \xi}}{\mathbb E^x[e^{\theta X_t}]}p_t(x,d\xi)
\end{equation}
with characteristic exponents
\[
\widetilde\phi(t,u)=\phi(t,u+\theta)-\phi(t,\theta),\quad \widetilde \psi(t,u)=\psi(t,u+\theta)-\psi(t,\theta).
\]
By differentiation at $t=0$ we obtain that its functional characteristics are given by
\[
\widetilde F(u)=F(u+\theta)-F(\theta)=F(u+\theta),\quad \widetilde R(u)=R(u+\theta)-R(\theta)=R(u+\theta).
\]

By construction, this process $(X,\mathbb Q_x)_{x\in D}$ is conservative, because $q_t(x,D)=1$, each $t,x$, and we have $\mathbb Q_x\sim \mathbb P_x$. Hence by ([AP3], Theorem 3.4), the only non-positive solution of \eqref{eq ric1} (with $\widetilde F, \widetilde R$
replacing $F,R$) for trivial initial data is the trivial one. Conversely, if we have two conservative processes $(X,\mathbb P_x)_{x\in D}$ and $(X,\mathbb Q_x)_{x\in D}$ with characteristics $(F,R)$ and $(\widetilde F,\widetilde R)$ and if
the process
\[
M_t^x=\exp\left(\theta(X_t-x)\right).
\]
is a local martingale, then $\mathbb Q_x\sim \mathbb P_x$ and $M_t$ is the associated density process, hence a martingale.

\subsubsection{Jump behaviour}
[AP7] reveals the surprising property that affine processes $X$ on $d\times d$ matrices of dimension $d>1$ do not exhibit jumps of infinite total variation. This finding is contrary to the situation on the positive real line. 
In the latter case, the linear jump intensities can be modulated down to $0$ at the boundary of the state-space, which is impossible in higher dimensions, due to the more complicated boundary structure of $S_d^+$, $d>1$.

Based on this finding, [AP7] proves that all jump-diffusions with non-degenerate diffusion component, are affine in the sense of \cite{dfs}, that is \eqref{eq: affine property} indeed holds for the characteristic function\footnote{The affine property in [AP3] is defined via the Laplace transform, which is the most convenient transform
for proper cones.}. For vanishing diffusion component, the processes are infinitely divisible, hence the same applies in this case. [AP9] extends this result to the Fourier-Laplace transform. It is still unknown whether for degenerate, non-zero diffusion coefficients $\alpha$, the characteristic function or the full Fourier-Laplace transform of affine processes on $S_d^+$ exhibits zeroes or not, unless
$X$ is a Wishart process with state-independent jumps (the matrix-variate generalization of the ``Basic Affine Jump Diffusions'' (BAJD) of \cite{duffie2001risk})

\subsubsection{Hypoellipticity}

In [AP8] we provide a new method to show absolute continuity of the transition function of affine jump-diffusions on $D=\mathbb R_+^m\times \mathbb R^n$.
 The existence of densities $g(t,x,\xi)$ of $X_t\mid X_0=x$ and their $C^k$ regularity allows us to expand them relative to a polynomial basis $p_\alpha$ \footnote{Here and below $\alpha$ is a multi-index, when the dimension of the state space is larger than one.} of a weighted $L^2$ space, say
\begin{equation}\label{proxy}
g(t,x,\xi)\approx \sum_{\vert \alpha \vert\leq N} c_\alpha p_\alpha.
\end{equation}
The basic idea to this stems from the polynomial character of affine jump-diffusions: Since the action of the infinitesimal generator of their Markovian semigroup reduces
to a linear operator acting on $\mathcal P_k$, one can determine finite moments of arbitrary order in rather explicit form. These moments, in turn,
determine the coefficients $c_\alpha$ of the  $L^2$ approximations.

A naive question in this regard would be why one performed density approximations, when one could actually get transition densities by mere Fourier Inversion of
the characteristic function of APs. There are, however, two main reasons why the use of our proxies do make sense.
\begin{itemize}
\item [(a)] Fourier-inversion, in particular for dimension $d\geq 2$, is numerically very challenging: For $d$ dimensions, one needs to solve a $d$-fold integral numerically.
\item [(b)] The (time-dependent) coefficients in the polynomial expansion of  $X_t\mid X_0=x$  can be calculated directly from a matrix exponential (of a constant matrix given in terms of the parameters of the affine process). Using the characteristic function for determining the density, however, forces
one to solve the generalized Riccati equations (a non-linear ODE as opposed to a linear IVP having an explicit solution given by the matrix exponential)
\end{itemize}

However, the existence of transition densities is a subtle issue itself. For solutions to SDEs of type
\[
dX_t=b(X_t)dt+\sigma(X_t)dB_t+dJ_t, \quad X_0\in\mathbb R^d,
\]
where $B$ is a $n$--dimensional standard Brownian motion, $b:\mathbb R^d\rightarrow\mathbb R^d$, $\sigma(x)$ is a $d\times m$ matrix for each $x$, and $J$  is a pure jump-process,
one can ask under which conditions $X_t$ ($t>0$) admits a Lebesgue density. Under certain regularity conditions, and in the case $J=0$,  \cite{Hoermander67} provides a sufficient and necessary one for the existence of densities. Writing $\sigma^0=b$, and $\sigma=(\sigma^1,\dots\sigma^n)$,
where $\sigma^i$ are column vectors, H{\"o}rmander's condition states that at any point $x\in D$, the Lie brackets
\[
\sigma^0,\quad  [\sigma^0,\sigma^1], \quad [[\sigma^0,\sigma^1],\sigma^2],\dots
\]
span $\mathbb R^d$. However, for affine jump-diffusions, $\sigma$ typically is only H\"older regular at the boundary of the state-space $D$, and derivatives lose regularity there as well.
Therefore,  H{\"o}rmander's  ingenious paper is not applicable for affine jump-diffusions, in general. 

In [AP8] we develop a Fourier method, which lets
us conclude regularity of the transition probabilities from the behaviour of the characteristic function
$\Phi(t,u,x)$ for $\|u\|\rightarrow\infty$. We apply the elementary fact
from Fourier Analysis \cite[Proposition 28.1]{Sato1999} which states, if
\begin{equation}\label{growth beh eq}
\lim_{\|u\|\rightarrow\infty}\frac{\|\Phi(t,u,x)\|}{\|u\|^k}=0,
\end{equation}
then $X_t\mid X_0=x$ admits a Lebesgue density $g$, and $g$ is $k$--times differentiable. To this end, we design a method, which
makes use of the particular form of the generalized Riccati Equations \eqref{eq ric1}--\eqref{eq ric2}. 

The sufficient conditions formulated in [AP8] to obtain growth behaviour \eqref{growth beh eq} are closely related to  H{\"o}rmander's conditions, when considered as translation into the Fourier domain. Current research concentrates on generalizations which allow for positive components with degenerate diffusion coefficients. In particular we have in mind transition probabilities for hazard rate models in default risk, which are naturally degenerate to some extent (e.g. the BAJD in \cite{duffie2001risk}).

\bibliographystyle{SIAM}
\bibliography{references}

\end{document}